\magnification=\magstep1
\baselineskip=14pt

\def\mapdown#1{\downarrow\rlap{$\vcenter{\hbox{$\scriptstyle#1$}}$}}

\centerline {\bf On the minimal free resolution of the universal ring}\centerline{\bf for
resolutions of length two.}
\bigskip
\centerline{ Andrew R. Kustin}
\centerline{Mathematics Department}
\centerline{University of South Carolina}
\centerline{Columbia, SC 29208}
\centerline{kustin@math.sc.edu}
\medskip
\centerline{and}
\medskip
\centerline{Jerzy M. Weyman\footnote*{Partially supported by NSF grant DMS-0300064.}}
\centerline{Mathematics Department}
\centerline{Northeastern University}
\centerline{ Boston, MA 02115}
\centerline{j.weyman@neu.edu}

\beginsection 0. Introduction.

Fix  positive integers $e$, $f$, and $g$, with $r_1\ge 1$ and $r_0\ge 0$, for $r_1$ and $r_0$ defined to be $f-e$ and $g-f+e$, respectively.  Hochster [Ho]  established the existence of a commutative noetherian ring $\tilde C$  and a universal resolution $$ U\!\!:\quad 0\to \tilde C^{e}\to \tilde C^{f}\to \tilde C^{g}$$ such that for any commutative noetherian ring $S$ and any resolution  $$ V\!\!: \quad0\to S^{e}\to S^{f}\to S^{g},$$ there exists a unique ring homomorphism $\tilde C\to S$ with $ V= U\otimes_{\tilde C} S$. Hochster showed that the universal ring $\tilde C$ is integrally closed and finitely generated as an  algebra over the ring of integers ${\bf Z}$. Huneke [Hu] identified the generators of $\tilde C$ as a ${\bf Z}$-algebra. These generators correspond to the entries of the two matrices from $ U$ and the ${g}\choose{r_1}$ multipliers from the factorization theorem  of Buchsbaum and Eisenbud [BE].
Bruns [B83] showed that $\tilde C$ is  factorial. Bruns [B84] also showed that  universal resolutions exist only for resolutions of length at most two. Heitmann [He] used Bruns' approach to universal resolutions in his counterexample to the rigidity conjecture.  
Pragacz and Weyman [PW] found the relations on the generators of $\tilde C$ and used Hodge algebra techniques to prove that ${\bf K}\otimes_{\bf Z}\tilde C$ has rational singularities when ${\bf K}$ is a field of characteristic zero.  
Tchernev  [T] used the theory of Gr\"obner bases to generalize and extend all of the above results with special interest in allowing an arbitrary  base ring $C_0$. In particular, his method yields the following results. 
\itemitem{\rm (a)} The ring $C_0\otimes_{\bf Z}\tilde C$ is factorial,  or Cohen-Macaulay, or  Gorenstein  if and only if the base ring $C_0$ has the same property.
\itemitem{\rm (b)} The ring  $C_0\otimes_{\bf Z}\tilde C$ is regular if and only if the base ring $C_0$ is regular and $r_1=1$.
\itemitem{\rm (c)} If $C_0$ is a perfect field of positive characteristic, then $C_0\otimes_{\bf Z}\tilde C$ is $F$-regular.

When $r_1$ is equal to 1, then the universal ring $\tilde C$ is the polynomial ring over $\bf Z$ with variables which represent entries of the second matrix from $ U$ together with variables which represent the Buchsbaum-Eisenbud multipliers. In particular, when $g=r_1=1$, then  the Hilbert-Burch theorem, which classifies all resolutions of the form $$0\to \tilde C^{f-1}\to \tilde C^{f}\to \tilde C^1,$$ is recovered. When $e=1$ and $r_0=0$, then the universal resolution looks like
$$0\to \tilde C^1\to \tilde C^{f}\to \tilde C^{f-1},$$
and
the universal ring $\tilde C$ is defined by the generic Herzog ideal of grade $f$   [AKM]  in the polynomial ring $\tilde B={\bf Z}[\{\phi_j\},\{\psi_{k,j}\}, a]$, where $1\le k\le f-1$ and $1\le j\le f$. The minimal resolution of $\tilde C$ by free $\tilde B$-modules is given in [KM]. 
  
The present paper concerns the universal ring $\tilde C$ when $r_0=0$. In this case,   $f=e+g$ and $\tilde C=\tilde B/\tilde  J$, for
$\tilde B$ equal to the polynomial ring ${\bf Z}[\{\phi_{j,i}\}, \{\psi_{k,j}\},a]$, with    $1\le i\le e$, $1\le j\le f$,
and $1\le k\le g$,
where  $\{\phi_{j,i}\}\cup\{\psi_{k,j}\}\cup \{a\}$  is a list of indeterminates over $\bf Z$.
The indeterminate $a$ corresponds to the unique Buchsbaum-Eisenbud multiplier which occurs in the present situation. Let
$\phi$ be the $f\times e$ matrix and $\psi$ be the $g\times f$ matrix  with
entries  $\phi_{j,i}$ and $\psi_{k,j}$, respectively. View the matrices $\phi$ and $\psi$ as  homomorphisms of $\tilde B$-modules:
 $$ \tilde B^e{\buildrel  \phi\over\rightarrow}  \tilde B^f {\buildrel  \psi\over\rightarrow}\tilde B^g.$$
 We give $\tilde  J$ in the language of [T].  For each $${\rm partition  
\ of\ \{1,\dots, f\}\ into\  I\cup \bar I\  with\  |I|=e\  and
\ |\bar I |=g,}\eqno{(*)}$$ let $\nabla_{\bar I,I}$ be the sign of the permutation which arranges the elements of $\bar I,I$ into increasing order, $\phi(I)$ the submatrix of $\phi$ consisting of the rows from $I$, and $\psi(\bar I)$ the submatrix of $\psi$ consisting of the columns from $\bar I$.   In this notation, the ideal $\tilde  J$ which defines the universal ring $\tilde C$ is   $$I_1(\psi\phi)+(\{\det \psi(\bar I)+\nabla_{\bar I,I} a\det \phi(I)\mid  I\cup \bar I\ {\rm from}\  ({*})\}).\eqno{(**)}$$  

One resolution, $\bf X$, of $\tilde C$ by free $\tilde B$-modules may be found in [K]. The resolution $\bf X$ is not minimal; but it is straightforward, coordinate free, and independent of characteristic. Furthermore, one can use $\bf X$ to calculate  ${\rm Tor}_{\bullet}^{\tilde B}(\tilde C,{\bf Z})$.
If $e$ and $g$ are both at least $5$, then ${\rm Tor}_{\bullet}^{\tilde B}(\tilde C,{\bf Z})$ is not a free abelian group; and therefore (see Roberts [R] or  Hashimoto [Ha]), the graded betti numbers in the minimal resolution of $\tilde C\otimes_{\bf Z} {\bf K}$ by free  $\tilde B\otimes_{\bf Z} {\bf K}$-modules depend on the characteristic of the field $\bf K$.

{\bf Henceforth, we work over a field $\bf K$ of characteristic zero.} Consider the vector spaces $E, F, G$ over $\bf K$ of dimensions $e, f,
g$ respectively, with $f=e+g$. 
Since we will apply the geometric technique of [W], we identify 
$B=\tilde B\otimes_{\bf Z}{\bf K}$ with the coordinate ring of the affine space 
$${\rm Hom}_{\bf K}(E,F)\times {\rm Hom}_{\bf K}(F,G)\times \bf K.$$The vector space ${\rm Hom}(E,F)$ is naturally equal to $F\otimes E^*$; and therefore,  $B$ is the polynomial ring 
$$B={\rm Sym}_{\bf K}(F^*\otimes_{\bf K} E)\otimes_{\bf K}{\rm Sym}_{\bf K}(G^*\otimes_{\bf K} F)\otimes_{\bf K}{\bf K}[a].$$
Let $E\otimes_{\bf K} B{\buildrel  \phi\over\rightarrow} F\otimes_{\bf K} B{\buildrel  \psi\over\rightarrow}  G\otimes_{\bf K} B$ be the natural  maps 
given by $$\phi(u)=\sum_i v_i\otimes (v_i^*\otimes u)$$ and $$\psi(v)=\sum_i w_i\otimes (w_i^*\otimes v),$$ for each $u\in E$ and $v\in F$. It is not necessary to pick bases; however, if $u_1,\dots,u_e$; $v_1,\dots,v_f$; and $w_1,\dots,w_g$ are bases  for the vector spaces $E$, $F$, and $G$; and    $u_1^*,\dots,u_e^*$; $v_1^*,\dots,v_f^*$; and $w_1^*,\dots,w_g^*$ are the corresponding dual bases for $E^*$, $F^*$, and $G^*$; then 
$\sum_{i} v_i\otimes v_i^*$, which is used in the definition of $\phi$, is the element in $F\otimes F^*$ which represents the identity map under the canonical identification of $F\otimes F^*$ with ${\rm Hom}(F,F)$. The coordinate functions in $B$ may be identified as 
$\phi_{i,j}=v^*_i\otimes u_j\in F^*\otimes E$ and $\psi_{i,j}=w^*_i\otimes v_j\in G^*\otimes F$. The matrices which represent the maps $\psi$ and $\phi$, with respect to the chosen bases, are the generic matrices   $(\psi_{i,j})$ and $(\phi_{i,j})$, respectively. We have $C=\tilde C\otimes_{\bf Z}{\bf K}$ and $J=\tilde JB$. So, $B$ is the polynomial ring ${\bf K}[\{\phi_{i,j}\}, \{\psi_{i,j}\}, a]$, $C=B/J$, and $J$ is given by (**).  In Theorem 5, we produce the minimal resolution $\bf G$ of $C$ by free $B$-modules. The ring $B$ is bigraded with $\phi_{i,j}\in B_{(1,0)}$, $\psi_{i,j}\in B_{(0,1)}$, and $a\in B_{(-e,g)}$. The ideal $J$  and the resolution $\bf G$  are homogeneous with respect to this bidegree. 

 Let us recall the result from [PW]
that gives a natural basis for the universal ring $C$.
We notice that the proper ${\rm GL}$-representation in $${C={{B}\over{J}}}=
{{{{\rm Sym}}_{\bf K}(F^*\otimes E)\otimes_{\bf K} {{\rm Sym}}_{\bf K}(G^*\otimes F)\otimes_{\bf K} {\bf K}[a]}\over{J}}$$ for the multiplier $a$ is $\bigwedge^eE^*\otimes \bigwedge^fF\otimes\bigwedge^gG^*$. Indeed, the representation   $S_{1^g}F\otimes S_{1^g}G^*$ is equal to 
$$ S_{\lambda}E\otimes S_{(\mu_1,\dots,\mu_{g-1},0,-\lambda_e,\dots,-\lambda_1)}F\otimes S_{\mu}G^*\otimes (\bigwedge^eE^*\otimes \bigwedge^fF\otimes\bigwedge^gG^*),$$for $\lambda=1^e$ and $\mu=0$. In other words, in $C$, each maximal minor of $\psi$ is equal to the appropriately signed complementary maximal minor of $\phi$ times the image of\hfil\break  $\bigwedge^eE^*\otimes \bigwedge^fF\otimes\bigwedge^gG^*$. 

\proclaim Remark. Notice that in {\rm[W]} one uses the notation $L_\lambda E$,
$K_\lambda E$ to denote the Schur and Weyl functors. In this paper we work
over a field of characteristic zero, so we have our $S_\lambda E$
isomorphic to $L_{\lambda^\prime}E$ or $K_\lambda E$, where
$\lambda^\prime$ is a conjugate partition. The module $S_\lambda E$ is defined for any dominant weight $\lambda$ {\rm(}i.e., for any integers $\lambda_1\ge \lambda_2\ge \dots\ge \lambda_e${\rm)} because $$S_{(\lambda_1,\dots,\lambda_e)}E=S_{(\lambda_1+t,\dots,\lambda_e+t)}E\otimes  (\bigwedge^eE^*)^{\otimes t}$$ for any integer $t$. 

\proclaim Proposition 1. The ring $C$ has the following decomposition to
representations of \hfil\break ${\rm GL}(E)\times {\rm GL}(F)\times {\rm GL}(G)${\rm :}
$$C=\oplus_{\lambda ,\mu, t}
 S_{\lambda}E\otimes S_{(\mu_1,\ldots ,\mu_{g-1},0,-\lambda_e ,\ldots ,-\lambda_1
)}F\otimes S_{\mu}G^*\otimes  (\bigwedge^e E^*\otimes
\bigwedge^f F \otimes
 \bigwedge^g G^*)^{\otimes t},$$
where we sum over all partitions $\lambda$ with $e$ parts, partitions $\mu$
with $g-1$ parts and $t\ge 0$.
Note that the representation corresponding to the triple $(\lambda ,\mu ,t)$
is a factor of \hfil\break $(S_\lambda E\otimes S_\lambda F^* )\otimes (S_\mu F\otimes
S_\mu G^* )\otimes a^t$.

\noindent{\sl Proof.} Applying Theorem 1.3 from [PW], or Theorem 5.8 from [T], we get

$$C=\oplus_{\lambda ,\mu, t}L_\lambda E\otimes L_{(e+g-\lambda_u  ,\ldots
,e+g-\lambda_1 , \mu_1 ,\ldots ,\mu_s )}F\otimes L_\mu
G^*\otimes(\bigwedge^g G^* )^{\otimes t}.$$
Changing Schur functors to Weyl functors (i.e., $L$'s to $S$'s), partitions $\lambda ,\mu$ to $\lambda' ,\mu'$ respectively, and adjusting powers of
determinant representations to get a ${\rm GL}(E)\times {\rm GL}(F)\times {\rm GL}(G)$ - equivariant
statement
 we get the
result. $\bullet$

\proclaim Corollary. The ring $C$ is a free ${\bf K}[a]$-module. 

The ring $C/aC$
is isomorphic to the factor of $A:={\bf K}[ \phi_{i,j}, \psi_{i,j}]$ by the ideal
$I$ given by the relations $\psi\phi=0$,
$\bigwedge^g \psi=0$. The ring $A=B/a$ inherits the bidegree of $B$ with  $\phi_{i,j}\in A_{(1,0)}$ and $\psi_{i,j}\in A_{(0,1)}$. 

In section one we recapitulate the geometric method for calculating syzygies. The resolution of $A/I$ by free $A$-modules is given in the second section. In section two, we also resolve a family of maximal Cohen-Macaulay modules over the determinantal ring $\bar A/I_g(\psi)$, for $\bar A={\rm Sym}_{\bullet}(F\otimes G^*)$. The familiar rank one maximal Cohen-Macaulay modules ${\rm Sym}_{i}({\rm cok}\, \psi)$, for $0\le i\le e+1$, which are resolved by the Eagon-Northcott complex, are members of our family.  Section three gives the resolution of the universal ring $C=B/J$ by free $B$-modules.   

\beginsection 1. Geometric technique of calculating syzygies.

\bigskip

In this section we provide a quick description of the main facts related to the 
geometric technique of calculating syzygies; see [W] for more details.
We work over a field $\bf K$. The characteristic of $\bf K$ must be zero for the Bott algorithm; otherwise, in this section, the characteristic of $\bf K$ is arbitrary. 

Let us consider the projective variety $V$ of dimension $m$.
Let $X=A^N_{\bf K}$ be the affine space. The space $X\times V$ can be
viewed
as a total space of trivial vector bundle $\cal E$ of dimension $N$ over
$V$.
Let us consider the subvariety $Z$ in $X\times V$ which is the total space
of a subbundle
$\cal S$ in $\cal E$.
We denote by $q$ the projection $q: X\times V\longrightarrow X$ and by
$q^\prime$ the restriction
of $q$ to $Z$. Let $Y=q(Z)$. We get the basic diagram
$$                  \matrix{ Z&\subset&X\times V\cr
                 \mapdown{q^\prime}&&\mapdown{q}\cr
                 Y&\subset&X\cr }$$

The projection from $X\times V$ onto $V$ is denoted by $p$ and the quotient
bundle ${\cal E}/{\cal S}$
by $\cal T$.
Thus we have the exact sequence of vector bundles on $V$
$$           0\longrightarrow {\cal S}\longrightarrow {\cal
E}\longrightarrow {\cal
T}\longrightarrow 0$$
 The dimensions of $\cal S$ and $\cal T$ will be denoted by $s$, $t$
respectively.
The coordinate ring of $X$ will be denoted by $A$. It is a polynomial ring
in
 $N$ variables over $\bf K$.
We will identify the sheaves on $X$ with $A$-modules.

 The locally free resolution of the sheaf ${\cal
O}_Z$ as an ${\cal O}_{X\times V}$-module is given by the Koszul complex
$$ {\cal K}_{\bullet} (\xi ): 0\rightarrow \bigwedge^t (p^* \xi
)\rightarrow \ldots
\rightarrow \bigwedge^2 (p^* \xi )\rightarrow p^* (\xi ) \rightarrow {\cal
O}_{X\times V}$$
where $\xi = {\cal T}^*$. The differentials in this complex are homogeneous
of degree $1$ in
the coordinate functions on $X$.
  The direct image $p_* ({\cal O}_Z )$ can be identified with the
the sheaf of algebras ${\rm Sym}(\eta )$ where $\eta = {\cal S}^*$.

The idea of the geometric technique is to use the Koszul complex ${\cal
K}(\xi )_{\bullet} $ to
construct for each vector bundle $\cal V$ on $V$ the free complex ${\bf
F}_{\bullet} ({\cal V})$ of
$A$-modules with the homology supported in $Y$. In many cases  the complex
${\bf F}({\cal O}_V )_{\bullet} $
gives the free resolution of the defining ideal of $Y$.

For every vector bundle $\cal V$ on $V$ we introduce the complex
$${\cal K}(\xi ,{\cal V})_{\bullet}  := {\cal K}(\xi )_{\bullet}
\otimes_{{\cal
O}_{X\times V}} p^* {\cal V}$$ This complex is a locally free resolution of
the ${\cal O}_{X\times V}$-module
$ M({\cal V}) :={\cal O}_Z \otimes p^* {\cal V}$.

Now we are ready to state the basic theorem (Theorem (5.1.2) in [W]).

\proclaim  Theorem 1. For a vector bundle $\cal V$ on $V$
we define a free graded $A$-module
$$ \ {\bf F}({\cal V})_i  = \bigoplus_{j\ge 0} H^j
(V,\bigwedge^{i+j}\xi\otimes{\cal V}
)\otimes_{\bf K} A(-i-j).$$
 \item {\rm (a)} There exist  minimal differentials
$$ d_i ({\cal V}): {\bf F}({\cal V})_i \rightarrow {\bf F}({\cal V})_{i-1}
$$
of degree $0$ such that ${\bf F}({\cal V})_{\bullet} $ is a complex of
graded free
$A$-modules with
$$ H_{-i} ({\bf F}({\cal V})_{\bullet} ) = {\cal R}^i q_* M({\cal V}).$$
In particular, the complex ${\bf F}({\cal V})_{\bullet} $ is exact in
positive degrees.
\item {\rm (b)} The sheaf ${\cal R}^i q_* M({\cal V})$ is equal to $H^i (Z,
M({\cal V}))$ and it
can be also identified with the graded $A$-module
$H^i (V, {\rm Sym} (\eta )\otimes {\cal V})$.
\item {\rm(c)} If $\phi : M({\cal V})\rightarrow M({\cal V}^\prime )(n)$ is a
morphism of
graded sheaves then there exists a morphism of complexes
$$f_{\bullet} (\phi ): {\bf F}({\cal V})_{\bullet} \rightarrow {\bf
F}({\cal V}^\prime )_{\bullet} (n).$$
Its induced map $H_{-i} (f_{\bullet} (\phi ))$ can be identified with the
induced map
$$H^i (Z, M({\cal V}) )\rightarrow H^i (Z, M({\cal V}^\prime ) )(n).$$

\vskip .2cm

If $\cal V$ is a one dimensional trivial bundle on $V$,  then the complex
${\bf F}({\cal V})_{\bullet} $
is denoted simply by ${\bf F}_{\bullet}$.

  The next theorem gives the criterion for the complex ${\bf F}_{\bullet}$
to be the free resolution of the
coordinate ring of $Y$.

\proclaim  Theorem 2. Let us assume that the map $q^\prime :
Z\longrightarrow Y$
is a birational isomorphism. Then the following properties hold.
\item {\rm (a)} The module $q^\prime_* {\cal O}_Z$ is the normalization of ${\bf
K}[Y]$.
\item {\rm (b)} If ${\cal R}^i q^\prime_* {\cal O}_Z = 0$ for $i>0$, then ${\bf
F}_{\bullet}$ is a finite
free resolution of the normalization of ${\bf K}[Y]$ treated as an
$A$-module.
\item {\rm (c)} If ${\cal R}^i q^\prime_* {\cal O}_Z = 0$ for $i>0$ and ${\bf
F}_0 = H^0 (V, \bigwedge^0 \xi
)\otimes A = A$, then $Y$ is normal and it has rational singularities.

This is Theorem (5.1.3) in [W].

The complexes $F({\cal V})_{\bullet} $ satisfy the Grothendieck type
duality. Let $\omega_V$ denote the
canonical divisor on $V$.

\proclaim Theorem 3. Let ${\cal V}$ be a vector bundle on $V$. Let us
introduce the dual bundle
$${\cal V}^\vee = \omega_V \otimes \bigwedge^t \xi^* \otimes {\cal V}^* .$$
Then
$$ F({\cal V}^\vee )_{\bullet}  = F({\cal V})_{\bullet}^* [m-t]$$

This is Theorem (5.1.4) in [W].

In all our applications the projective variety $V$ will be a Grassmannian.
To fix the notation, let us work with the
Grassmannian ${\rm Grass}(r, E)$ of subspaces of dimension $r$ in a vector space
$E$ of dimension $n$.
Let $$0\rightarrow {\cal R}\rightarrow E\times {\rm Grass}(r, E)\rightarrow {\cal
Q}\rightarrow 0$$
be a tautological sequence of the vector bundles on ${\rm Grass}(r, E)$.

Assume that the characteristic of the field $\bf K$ is zero.
Then the vector bundle
$\xi$ will be a direct sum of the bundles of the form $S_{\lambda_1 ,\ldots
,\lambda_{n-r}}{\cal Q}\otimes S_{\mu_1
,\ldots ,\mu_r}{\cal R}$.  Thus all the exterior powers of $\xi$ will also
be the direct sums of such bundles. We will
apply repeatedly the following result to calculate cohomology of vector
bundles $S_{\lambda_1 ,\ldots
,\lambda_{n-r}}{\cal Q}\otimes S_{\mu_1 ,\ldots ,\mu_r}{\cal R}$.

\proclaim Proposition 2 (Bott's algorithm). Assume that the characteristic
of $\bf K$ is zero.
The cohomology of the vector bundle $S_{\lambda_1 ,\ldots
,\lambda_{n-r}}{\cal Q}\otimes
S_{\mu_1 ,\ldots ,\mu_r}{\cal R}$
 on ${\rm Grass}(r, E)$ is calculated
as follows. We look at the weight
$$(\lambda ,\mu ) = (\lambda_1 ,\ldots ,\lambda_{n-r},\mu_1 ,\ldots ,\mu_r
) $$
 and add to it $\rho = (n, n-1,\ldots
,1)$.  If the resulting sequence
$$(\lambda ,\mu ) +\rho =(\lambda_1 +n,\ldots ,\lambda_{n-r} +r+1,\mu_1 +r
,\ldots
,\mu_r +1)$$
 has repetitions, then
$$ H^i ({\rm Grass}(r, E), S_\lambda {\cal Q}\otimes S_{\mu} {\cal R} )=0$$
for all $i\ge 0$.
If the resulting sequence has no repetitions, there is a unique permutation
$w\in \Sigma_n$
that makes this sequence decreasing.  Then the sequence $\nu =w((\lambda
,\mu ) +\rho )-\rho$ is again a non-increasing
sequence. Then the sheaf $S_{\lambda}{\cal Q}\otimes S_{\mu} {\cal R}$ has
only one non-zero
cohomology group, the group $H^l$ where $l=l(w)$ is the length of $w$.
This cohomology  group is isomorphic to the representation
$S_\nu E$ of ${\rm GL}(E)$ corresponding
to the highest weight $\nu$ {\rm(}the so-called Schur module{\rm)}.

This is Corollary (4.1.9) in [W].

\beginsection 2. The resolution of $A/I$.

We apply the geometric technique to calculate the minimal free resolution
of $A/I$ as an $A$-module. The notation is set up in the final Corollary of the Introduction. Recall that $\bf K$ is a field of characteristic zero. We use freely the notation of [W].
Denote $$X=\lbrace (d_2 ,d_1 )\in {\rm Hom}_{\bf K} (E,F)\times {\rm Hom}_{\bf K} (F,G)\rbrace .$$
Therefore we have $A={\bf K}[X]$.
Consider the incidence variety
$$Z=\lbrace (d_2 ,d_1 , {R} )\in X\times {\rm Grass}(e+1,F)\ |\
{\rm Im}(d_2)\subseteq {R}\subseteq {\rm Ker}(d_1)\ \rbrace .$$
Clearly the image $q(Z)$ by the first projection $q:Z\rightarrow X$ is
equal to the set $Y:=V(I)$.
Notice that $Z$ is the desingularization of $Y$ because generically on $Y$
we have ${R}={\rm Ker}(d_1 )$ and the projection $q$ is obviously proper.

We are in the situation described in the previous section. In this special
case we have $\xi = E\otimes {\cal Q}^*\oplus {\cal R}\otimes G^* $. We
also have $\eta =E\otimes {\cal R}^*\oplus {\cal Q}\otimes G^*$.
Let us look at the cohomology groups of the exterior powers of $\xi$ and of
symmetric powers of $\eta$.

\proclaim Proposition 3.  We have
\item{\rm (a)} $H^i ({\rm Grass} (e+1,F), {\rm Sym}_j (\eta ))=0$ for $i>0$,
\item{\rm (b)} $H^0 ({\rm Grass}(e+1,F), {\rm Sym}_j (\eta ))= (A/I)_j$ for all $j\ge 0$.

\noindent{\sl Proof.} We have
$${\rm Sym} (\eta )=\oplus_{\lambda ,\mu}S_\lambda E\otimes S_\lambda {\cal
R}^*\otimes S_\mu {\cal Q}\otimes S_\mu G^*$$
where we sum over partitions $\lambda$ with $e$ parts and partitions $\mu$
with $g-1$ parts. We notice that higher cohomology of the bundles
$ S_\lambda {\cal R}^*\otimes S_\mu {\cal Q}$ is zero, with $H^0$ being
just $S_{(\mu_1 ,\ldots ,\mu_{g-1} ,0 ,-\lambda_{e},\ldots ,-\lambda_1 )}
F$. Comparing it with Proposition 1 we are done.$\bullet$

Proposition 3 implies that the complex ${\bf F}_\bullet$ is a minimal free
resolution of the coordinate ring of $Y$.

Let us analyze the cohomology of the exterior powers of $\xi$. We have

$$\bigwedge^\bullet (\xi )=\oplus_{\lambda ,\mu}\ S_{\lambda' }E\otimes
S_\lambda {\cal Q}^*\otimes S_\mu {\cal R}\otimes S_{\mu' }G^* .$$

To calculate the cohomology of the summand corresponding to the pair
$(\lambda ,\mu )$ we need to apply the Bott algorithm to the sequence
$$(-\lambda_{g-1},\ldots ,-\lambda_1 ,\mu_1 ,\ldots ,\mu_{e+1}) .$$

\proclaim Proposition 4.
\item{\rm (a)} The representations of $F$ we get from the above procedure are
all of the type $\bigwedge^s F$ ($0\le s\le f)$.
\item{\rm (b)} The ring ${\bf K}[Y]$ is normal and  Gorenstein, of codimension $eg+1$
in $X$.

\noindent{\sl Proof.} Let us look what will be the highest number in our
sequence after applying  Bott's algorithm. It clearly is either
$-\lambda_{g-1}$ or $\mu_1 -g+1$. But $\mu_1\le g$, otherwise the
corresponding summand is zero as it involves the factor $S_{\mu'}G^*$. Thus
the first number is $\le 1$. Similarly, the last number is either
$\mu_{e+1}$ or $-\lambda_1 +e+1$. Since $\lambda_1\le e$ (otherwise the
summand is zero, as it contains factor $S_{\lambda'}E$), we see that the
last number is $\ge 0$.
Thus our weight has to be of the type $(1^s ,0^{f-s})$.

Let us look at the top exterior power of $\xi$. Clearly this is
$$\bigwedge^{top}\xi =S_{(g-1)^e}E\otimes S_{e^{g-1}}{\cal Q}^*\otimes
S_{g^{e+1}}{\cal R}\otimes S_{(e+1)^g}G^* .$$
To calculate the corresponding term, we need to apply Bott's algorithm to
the sequence
$(-e^{g-1}, g^{e+1})$ which gives the representation $\bigwedge^f F$ in
$H^{(g-1)(e+1)}$. This is the top of the resolution. The representation
there is
$$\bigwedge^{top}\xi =S_{(g-1)^e}E\otimes \bigwedge^f F\otimes
S_{(e+1)^g}G^* $$
in the homological degree $e(g-1)+g(e+1)-(g-1)(e+1)=eg+1$. The
representation is one dimensional, therefore ${\bf K}[Y]$ is Gorenstein, of
codimension $eg+1$ as claimed. The normality follows because $Z$ is a
desingularization of $Y$.$\bullet$

\proclaim Remarks.
\item{\rm (a)} The strand of the complex ${\bf F}_\bullet$ with the
${\rm SL}(F)$-component $\bigwedge^0 F$ is just
$$\oplus_\lambda S_{\lambda' }E\otimes S_{\lambda}G^* .$$
This is just a subcomplex of the Koszul complex on the composition
$\psi\phi$.
\item{\rm (b)} The differential in the complex ${\bf F}_\bullet$ has three
components. One involves only the map $\phi$, the other only the map
$\psi$, and the third component
is of degree $(1,1)$ in $\phi$ and $\psi$, and it does not change the
$F$-component of the term. We refer to these components as respectively
$\phi$-component, $\psi$-component and $(\psi\phi )$-component.

 Let us look at the terms in the complex ${\bf F}_\bullet$. One way to do
that is to look at the terms with a fixed $\lambda$.
In order to describe this part of the complex we need another geometric
construction related to the Grassmannian of $G$.
Consider ${\rm Grass}(g-1, G)$ with the tautological sequence
$$0\rightarrow {\bar{\cal R}}\rightarrow G\times {\rm Grass}(g-1, G)\rightarrow
{\bar{\cal Q}}\rightarrow 0.$$
We are dealing with the polynomial ring ${\bar A}={\rm Sym} (F\otimes G^* )$ and
the modules supported in the determinantal varieties of maps $\psi$ of rank
$\le g-1$. We look at twisted complexes ${\bar{\bf F}}(S_\nu {\bar{\cal
R}}^* )_\bullet$ which come from taking $\xi =F\otimes {\bar{\cal Q}}^*$.
Each such complex is the pushdown the locally free resolution of the sheaf
$${\cal M}(\nu ):= S_\nu {\bar{\cal R}}^*\otimes {\rm Sym} (F\otimes {\bar{\cal
R}}^* ) .$$

\proclaim Proposition 5. The sheaf ${\cal M}(\nu )$ has no higher
cohomology. Thus the complex ${\bar{\bf F}}(\nu )_\bullet$ is a free
resolution of the $\bar A$-module
$$M(\nu ):= H^0 ({\rm Grass}(g-1 ,G), {\cal M}(\nu ) ).$$
Assume that $\nu\subset e^{g-1}$. Then the complex ${\bar{\bf F}}(\nu
)_\bullet$ is a complex of length $f-g+1$. Thus the corresponding module
$M(\nu )$ is a maximal Cohen-Macaulay module.

\noindent{\sl Proof.} This is a standard application of geometric
technique, see [W], ch. 6.$\bullet$

\proclaim Remark. Let us look at the resolution of $M(\nu )$ more
precisely. It is a pushdown of the twisted Koszul complex
$$S_\nu {\bar{\cal R}}^*\otimes \bigwedge^\bullet (F\otimes {\bar{\cal
Q}}^* ) .$$
Thus we can describe  the terms as $\bigwedge^i F$ tensored with the
representation $S_{\mu (i)}G$ where $\mu (i)$ is the result of Bott
algorithm applied to the weight
$$(-i, -\nu_{g-1},\ldots ,-\nu_1 ).$$
The terms we get in $H^0$ correspond to $i$ satisfying $-i\ge -\nu_{g-1}$. For each such $i$, the $H^0$-module is equal to 
$$\bigwedge^i F\otimes S_{(-i, -\nu_{g-1},\ldots ,-\nu_1 )}G=\bigwedge^i F\otimes S_{(\nu_1,\ldots,\nu_{g-1},i)}G^*,$$and it 
 appears in the $i$-th place in the complex ${\bar{\bf F}}(\nu
)_\bullet$.
The terms we get in $H^s$ for $s\ge 1$  correspond to $i$ satisfying the
inequalities
$$-\nu_{g-s}-1\ge -i+s\ge-\nu_{g-s-1}.$$ For each pair $(i,s)$, the $H^s$-module is equal to 
$$\matrix{\bigwedge^i F\otimes S_{(-\nu_{g-1}-1,\ldots , -\nu_{g-s}-1, -i+s,
-\nu_{g-s-1},\ldots ,-\nu_1 )}G\cr =\bigwedge^i F\otimes S_{(\nu_1,\ldots ,\nu_{g-s-1},i-s,\nu_{g-s}+1,\ldots ,\nu_{g-1}+1)}G^*,}$$ and it appears in the $i-s$'th place in the complex ${\bar{\bf F}}(\nu
)_\bullet$.

\proclaim  Proposition 6. Let $\nu$ be a partition contained in the
rectangle $e^{g-1}$. Then the terms of the complex ${\bf F}_\bullet$
containing the factor $S_{\nu'} E$ are identical with the terms of the
complex $S_{\nu'} E\otimes {\bar{\bf F}}(\nu)_\bullet [|\nu |]$. Here $[i]$
means homological shift, i.e., the  term in position zero of
 $S_{\nu'} E\otimes {\bar{\bf F}}(\nu)_\bullet $ occurs in ${\bf F}_{|\nu
|}$.

\noindent{\sl Proof.} Direct calculation. The lowest term where $S_\lambda
E$ occurs was described in the remark (a) preceding  Proposition
5.$\bullet$

The modules $M(\nu)$ of Proposition 5 are maximal Cohen-Macaulay modules over the determinantal ring $\bar A/I_g(\psi)$, 
for $\bar A= {\rm Sym}(F\otimes G^*)$, where $\psi : F\otimes_{\bf K} \bar A\to G\otimes_{\bf K} \bar A$ is the natural map. These modules have independent interest. In Theorem 4 we record the $\bar A$ resolution ${\bf t}_{\nu}$ of $H_0({\bf t}_{\nu})=M(\nu)$ using one parameter $k$ in place of the two parameters $i$ and $s$ that were used to date. Recall that $\bf K$ is a field of characteristic zero, $F$ and $G$ are vector spaces over $\bf K$ of dimension $f$ and $g$, respectively, and $e=f-g$. 

\proclaim Definition. Let $k$ be an integer and $\nu=(\nu_1,\dots,\nu_{g-1})$ be a dominant weight. 

\noindent{\rm(a)} Let $i=\nu_k'$, which is defined to be the number of indices $j$ with $\nu_j\ge k$.  Notice that $\nu_i\ge k> \nu_{i+1}$.  Define $p(\nu;k)$ to be the dominant weight  
$$p(\nu;k)=(\nu_1,\dots,\nu_{i},k,\nu_{i+1}+1,\dots,\nu_{g-1}+1),$$
$N(\nu;k)$ to be the integer $g-1-\nu_k'+k$, and 
 and $t_{\nu;k}$ to be the free $\bar A$-module 
$$t_{\nu;k}=\bigwedge^{N(\nu;k)}F\otimes_{\bf K}  S_{p(\nu;k)}G^*\otimes_{\bf K} \bar A.$$

\noindent{\rm (b)} Define a homomorphism $t_{\nu;k}\to t_{\nu;{k-1}}$.
Let $N=1+\nu_{k-1}'-\nu_k'$. It follows that there exist dominant weights $\alpha$ and $\beta$ with $\alpha_{\rm{last}}\ge k> \beta_1$,
$$p(\nu;k)=(\alpha,k^N,\beta),\quad{\rm {and}}\quad p(\nu;k-1)=(\alpha,(k-1)^N,\beta).$$
The homomorphism $t_{\nu;k}\to t_{\nu;{k-1}}$ is  the composition

$$t_{(\nu;k)}=\bigwedge^{N(\nu;k)}F\otimes S_{p(\nu;k)}G^*\otimes \bar A
\to S_{1^{N}}G^* \otimes S_{p(\nu;k-1)}G^*\otimes \bar A\to $$$$
\bigwedge^{N(\nu;k)}F\otimes S_{1^{N}} F^*\otimes S_{p(\nu;k-1)}G^*\otimes \bar A\to 
 \bigwedge^{N(\nu;k)-N}F\otimes  S_{p(\nu;k-1)}G^*\otimes \bar A=t_{(\nu;k-1)},$$where the first map is the Pieri map, the second is $\bigwedge^N\psi^*$, and the third is the module action of $\bigwedge^{\bullet}F^*$ on $\bigwedge^{\bullet}F$.

\noindent{\rm (c)} For each dominant weight  $\nu=(\nu_1,\dots,\nu_{g-1})$, we define
the complex ${\bf t}_{\nu}$:
$$ \dots\to t_{\nu;k}\to t_{\nu;{k-1}}\to \dots,$$
with $t_{\nu;k}$ in position $k$.

\proclaim Remarks. 

\noindent{\rm  (a)} The dominant weight  
$p(\nu;k)$ may be interpreted as the result of applying Bott's algorithm to the sequence
$$\nu_1,\dots,\nu_{g-1},N(\nu;k).$$

\noindent {\rm  (b)} If $\nu_{g-1}\ge -1$ and $k<0$, then  $t_{\nu;k}=0$.  
 
\noindent {\rm  (c)} If $\nu_{g-1}\ge -1$ and $k\ge 0$, 
then $p(\nu;k)$ is a partition.

\noindent {\rm  (d)} The maps and modules of ${\bf t}_{\nu}$ form a complex because the Littlewood-Richardson rule tells us that the only coordinate free $\bf K$-vector space map 
$$S_{\alpha,k^{1+N},(k-1)^M,\beta}G^*\to S_{\alpha,(k-1)^{N},(k-2)^{1+M},\beta}G^*\otimes S_{1^{2+N+M}}G^*$$is zero, when $\alpha$ and $\beta$ are dominant weights with $\alpha_{\rm{last}}\ge k$ and $k-1>\beta_1$.

\proclaim Observation. If $\nu=(\nu_1,\dots,\nu_{g-1})$ is a dominant weight  and $\mu=(e-\nu_{g-1},\dots,e-\nu_1)$, then the complexes ${\bf t}_{\nu}$ and $({\bf t}_{\mu})^*[-e-1]$ are isomorphic. Furthermore, if $\nu\subset [-1,e+1]^{g-1}$, then 
$\mu$ also sits in $[-1,e+1]^{g-1}$ and $({\bf t}_{\nu})_i=0$ for $i<0$ or $e+1<i$.   

\noindent{\sl Proof.}  Let $k$ and $\ell$ be integers with $k+\ell=e+1$. The modules 
$$ \bigwedge^{N(\nu;k)}F\otimes S_{p(\nu;k)}G^*\quad{\rm {and}}\quad \bigwedge^{N(\mu;\ell)}F\otimes S_{p(\mu;\ell)}G^*$$
 are dual to one another because $N(\nu;k)+N(\mu;\ell)=f$ and if 
$$p(\nu;k)=(A_1,\dots,A_{g})\quad{\rm {and}}\quad p(\mu;\ell)=(B_1,\dots,B_g),$$then $A_i+B_{g+1-i}=e+1$. 
 All of the other assertions may be readily verified. $\bullet$ 

\proclaim  Theorem 4. If  $\nu_{g-1}\ge -1$, then 
\item{\rm (a)}${\bf t}_{\nu}$ is a resolution of $H_0({\bf t}_{\nu})$, and
\item{\rm (b)}$H_0({\bf t}_{\nu})$ is a module over $\bar A/I_g(\psi)$.

\proclaim  Corollary. If  $\nu\subset [-1,e+1]^{g-1}$, then
\item{\rm (a)}$H_0({\bf t}_{\nu})$ is a perfect $\bar A$-module with 
$${\rm Ext}_{\bar A}^{e+1}(H_0({\bf t}_{\nu}),\bar A)=H_0({\bf t}_{\mu})$$ for $\mu=(e-\nu_{g-1},\dots,e-\nu_1)$, and
\item{\rm (b)}$H_0({\bf t}_{\nu})$ is a maximal Cohen-Macaulay  $\bar A/I_g(\psi)$-module.

\proclaim Example. In particular, if $\nu=(i^{g-1})$, then the complex ${\bf t}_{\nu}$ is isomorphic to the Eagon-Northcott complex ${\cal C}^i$, see, for example, {\rm [E, Figure A2.6]}, and $$H_0({\bf t}_{i^{g-1}})= \cases {\bigwedge^{e+1}{\rm cok}(\psi^*)&{\rm if}  $i=-1$\cr
\bar A/I_g(\psi)&{\rm if} $i=0$\cr
{\rm Sym}_i({\rm cok}(\psi))&{\rm if} $1\le i$.}$$

We return to the resolution ${\bf F}_{\bullet}$. We incorporate the idea in Proposition 6 together with the notation of Theorem 4.  For each partition $\nu$ and each integer $k$, let 
$$T_{\nu;k}= S_{\nu'}E\otimes_{\bf K} \bigwedge^{N(\nu;k)}F\otimes_{\bf K} S_{p(\nu;k)}G^*\otimes_{\bf K}A.$$
We have established the following result. 

\proclaim Proposition 7. 
The terms of the complex ${\bf F}_\bullet$ are the
following
$$ {\bf F}_\bullet =\oplus_{\nu ;k} T_{\nu;k} (-|\nu|,-|\nu|-N(\nu;k)),$$
where we sum over all partitions $\nu$ contained in $e^{g-1}$ and integers $0\le k\le e+1$. The term $T_{\nu:k} (-|\nu|,-|\nu|-N(\nu;k))$ appears in ${\bf F}_{|\nu |+k}$.

\proclaim Example. Let us take $e=g=2$. We give two versions of our resolution ${\bf F}_\bullet$. In the first version, our resolution has the following
terms where we write
$(a,b;c;d,e)$ for $S_{(a,b)}E\otimes \bigwedge^cF\otimes
S_{(d,e)}G^*$.
$$\matrix{(1,1;4;3,3)\otimes A(-2,-6)\cr
\downarrow \cr
(1,0;4;3,2)\otimes A(-1,-5)\oplus (1,1;2;2,2)\otimes A(-2,-4)\cr
\downarrow\cr
(1,0;3;2,2)\otimes A(-1,-4)\oplus (1,1;1;2,1)\otimes A(-2,-3)\oplus (0,0;4;3,1)\otimes A(0,-4)\cr
\downarrow\cr
(1,1;0;2,0)\otimes A(-2,-2)\oplus (0,0;3;2,1)\otimes A(0,-3)\oplus (1,0;1;1,1)\otimes A(-1,-2)\cr
\downarrow\cr
(0,0;2;1,1)\otimes A(0,-2)\oplus (1,0;0;1,0)\otimes A(-1,-1)\cr
\downarrow\cr
(0,0;0;0,0)\otimes A}$$

The terms of ${\bf F}_\bullet$ are also listed in the following picture, which has the added advantage of giving insight into the maps of ${\bf F}_\bullet$. The row which corresponds to the partition $\nu$ is $S_{\nu'}E\otimes_{\bf K} {\bf t}_{\nu}$ as described in Theorem 4. Each row is acyclic. The Koszul complex map down the column on the right induces an acyclic sequence on the zeroth homology of the rows. An iterated mapping cone produces the complex ${\bf F}_\bullet$. In other words, there is a map of complexes from the middle row to the bottom row; there is a map of complexes from the top row to the mapping cone formed from the bottom two rows; and  ${\bf F}_\bullet$ is the mapping cone of this second map of complexes.  Notice that it is not correct to think of this picture as a double complex. The ``knight move'' $T_{2;1}(-2,-3)\to T_{0;2}(0,-3)$ which is induced by $\bigwedge^2 \phi$, is one of the components of the differential of ${\bf F}_\bullet$.
$$ 
\matrix{
T_{2;3}(-2,-6)&\to&T_{2;2}(-2,-4)&\to&T_{2;1}(-2,-3)&\to&T_{2;0}(-2,-2)\cr
&&&&&&\mapdown{}\cr
T_{1;3}(-1,-5)&\to&T_{1;2}(-1,-4)&\to&T_{1;1}(-1,-2)&\to&T_{1;0}(-1,-1)\cr
&&&&&&\mapdown{}\cr
T_{0;3}(0,-4)&\to&T_{0;2}(0,-3)&\to&T_{0;1}(0,-2)&\to&T_{0;0}}$$

\proclaim Example. Let us take $e=2, g=3$. Our resolution has the following
terms where we write
$(a,b;c;d,e,f)$ for $S_{(a,b)}E\otimes \bigwedge^cF\otimes
S_{(d,e,f)}G^*$.

$$\matrix{(2,2;5;3,3,3)\otimes A(-4,-9)\cr
\downarrow \cr
(2,1;5;3,3,2)\otimes A(-3,-8)\oplus (2,2;2;2,2,2)\otimes A(-4,-6)\cr
\downarrow\cr
(2,2;1;2,2,1)\otimes A(-4,-5)\oplus (2,1;3;2,2,2)\otimes A(-3,-6)\oplus  (2,0;5;3,2,2)\otimes A(-2,-7)\cr \oplus (1,1;5;3,3,1)\otimes A(-2,-7)\cr
\downarrow\cr
(2,2;0;2,2,0)\otimes A(-4,-4)\oplus (2,0;4;2,2,2)\otimes A(-2,-6)\oplus
(1,0;5;3,2,1)\otimes A(-1,-6)\cr \oplus (1,1;3;2,2,1)\otimes A(-2,-5)\oplus (2,1;1;2,1,1)\otimes A(-3,-4)\cr
\downarrow\cr
(1,0;4;2,2,1)\otimes A(-1,-5)\oplus (1,1;2;2,1,1)\otimes A(-2,-4)\oplus
(2,1;0;2,1,0)\otimes A(-3,-3)\cr \oplus (0,0;5;3,1,1)\otimes A(0,-5)\oplus (2,0;1;1,1,1)\otimes A(-2,-3)\cr
\downarrow\cr
(0,0;4;2,1,1)\otimes A(0,-4)\oplus (1,0;2;1,1,1)\otimes A(-1,-3)
\oplus (1,1;0;2,0,0)\otimes A(-2,-2)\cr \oplus (2,0;0;1,1,0)\otimes A(-2,-2)\cr
\downarrow\cr
(0,0;3;1,1,1)\otimes A(0,-3)\oplus (1,0;0;1,0,0)\otimes A(-1,-1)\cr
\downarrow\cr
(0,0;0;0,0,0)\otimes A}$$

\filbreak
Also, ${\bf F}_\bullet$ is the iterated mapping cone of a picture built using the following modules. 
$$ 
\matrix  
{T_{2,2;3}(-4,-9)&\to&T_{2,2;2}(-4,-6)&\to&T_{2,2;1}(-4,-5)&\to&
T_{2,2;0}(-4,-4)\cr
&&&&&&\mapdown{}\cr
T_{2,1;3}(-3,-8)&\to&T_{2,1;2}(-3,-6)&\to&T_{2,1;1}(-3,-4) &\to&T_{2,1;0}(-3,-3)\cr
&&&&&&\mapdown{}\cr
\matrix{ T_{1,1;3}(-2,-7)\cr\oplus\cr T_{2,0;3}(-2,-7}&\to&\matrix{T_{1,1;2}(-2,-6) \cr\oplus\cr T_{2,0;2}(-2,-5}&\to&
\matrix{T_{1,1;1}(-2,-3)\cr\oplus\cr T_{2,0;1}(-2,-4}
&\to& \matrix {T_{1,1;0}(-2,-2)\cr\oplus\cr T_{2,0;0}(-2,-2}\cr
&&&&&&\mapdown{}\cr
T_{1,0;3}(-1,-6)&\to&T_{1,0;2}(-1,-5)&\to&T_{1,0;1}(-1,-3) &\to&T_{1,0;0}(-1,-1)\cr
&&&&&&\mapdown{}\cr
T_{0,0;3}(0,-5)&\to&T_{0,0;2}(0,-4)&\to&T_{0,0;1}(0,-3)&\to&T_{0,0;0}}$$

\proclaim Remarks.
\item{\rm (a)} The $\psi$ component of the complex ${\bf F}_\bullet$ is just the
sum of differentials in in the complexes ${\bar{\bf F}}(\nu )_\bullet$.
\item{\rm (b)} The complementary partitions with respect to the rectangle
$e^{g-1}$ give the parts ${\bar{\bf F}}(\nu )_\bullet$ that are dual to
each other. Thus we can see in this way that ${\bf F}_\bullet$ is
self-dual. Also it length is $\le e(g-1)+f-g+1=e(g-1)+e+1=eg+1$, and thus
it is equal to $eg+1$.

The description of the terms of the complex ${\bf F}_\bullet$ given in
Proposition 6 is not accidental. It comes from a pushdown of different
Koszul complex.
Consider still the Grassmannian ${\rm Grass}(g-1 ,G)$ with the tautological
sequence
$$0\rightarrow {\bar{\cal R}}\rightarrow G\times {\rm Grass}(g-1, G)\rightarrow
{\bar{\cal Q}}\rightarrow 0.$$
Consider the sheaf of algebras
$${\cal B}={\rm Sym} (E\otimes F^* )\otimes {\rm Sym} (F\otimes {\cal R}^* )$$
over ${\rm Grass}(g-1 ,G)$. Obviously we have linear maps
$$\phi :E\otimes {\cal B}\rightarrow F\otimes {\cal B},\quad  \psi' :F\otimes
{\cal B}\rightarrow {\cal R}\otimes{\cal B}.$$
of sheaves of $\cal B$-modules. Consider now the condition $\psi'\phi =0$.
The Koszul complex given by the entries of the composition is acyclic. Thus
we get an acyclic complex of sheaves of $\cal B$-modules
$${\cal K}_\bullet :0\rightarrow {\cal K}_{e(g-1)}\rightarrow {\cal
K}_{e(g-1)-1}\rightarrow\ldots\rightarrow {\cal K}_1\rightarrow {\cal
K}_0$$
with ${\cal K}_i =\bigwedge^i (E\otimes {\cal R}^* )\otimes {\cal B}$.
Notice that
$${\cal K}_i=\oplus_{|\nu |=i} S_{\nu'}E\otimes S_\nu {\cal R}^*\otimes
{\cal B}.$$
Let us denote ${\hat{\cal M}}(\nu ):= S_\nu {\cal R}^*\otimes {\cal B}$ and
${\hat M}(\nu ):= H^0 ({\rm Grass}(g-1, G), {\hat{\cal M}}(\nu ))$.

\proclaim Proposition 8. We have the following properties
\item{\rm(a)} $H^j ({\rm Grass}(g-1, G), {\cal K}_i )=0$ for $j>0$, $0\le i\le
e(g-1)$, $H^j ({\rm Grass}(g-1, G), {\hat{\cal M}}(\nu ))=0$, for $j>0$.
\item{\rm(b)} The resolution of ${\hat M}(\nu )$ as an $A$-module is just
${\bar{\bf F}}(\nu )_\bullet\otimes_{\bar A}A$.

\noindent{\sl Proof.} This is clear from the definitions. $\bullet$

The Koszul complex ${\cal K}_\bullet$ induces an acyclic complex of
sections
$$0\rightarrow K_{e(g-1)}\rightarrow
K_{e(g-1)-1}\rightarrow\ldots\rightarrow K_1\rightarrow K_0 .$$
where $K_i := H^0 ({\rm Grass}(g-1, G), {\cal K}_i )=\oplus_{|\nu |=i}
S_{\nu'}E\otimes {\hat M}(\nu )$. We can now use the iterated mapping cone
construction to construct the resolution ${\bf F}'_\bullet$ of the zero-th
homology group of $K_\bullet$. The terms of this resolution are the same as
the terms of ${\bf F}_\bullet$. The whole process can be made ${\rm GL}(E)\times
{\rm GL}(F)\times {\rm GL}(G)$-equivariant.

\proclaim Proposition 9. The resulting complex ${\bf F}'_\bullet$ is
isomorphic to ${\bf F}_\bullet$.

\noindent{\sl Proof.} Both complexes have the same terms and are
${\rm GL}(E)\times {\rm GL}(F)\times {\rm GL}(G)$-equivariant and thus minimal (every
representation occurs in ${\bf F}_\bullet$ at most once). Looking at the
first and zero-th term of ${\bf F}'_\bullet$ we can identify the
differential as the $g\times g$ minors of $\psi$ and the entries of the
composition matrix $\psi\phi$.$\bullet$

The description of the terms of the complex ${\bf F}_\bullet$ given in
Proposition 7 allows us also to understand the $\phi$ component of the
differential.
Consider two terms of ${\bf F}_\bullet$ with the same factor $S_\mu G^*$,
but occurring in neighboring cohomology groups. In other words, we are given the data $\nu;k$ from Proposition 7. Let $j=\nu_k'$. Assume $k\ge 1$ and $j\ge 1$. Let $\rho$ equal $\nu$ with $\nu_j$ replaced by $k-1$. One may check that $\rho$ is a partition, $p(\nu;k)=p(\rho,\nu_j)$, 
$\rho_{\nu_j}'=j-1$, and $N(\rho;\nu_j)=N(\nu;k)+\nu_j-k+1$. The map
$$T_{\nu;k}\to T_{\rho;\nu_j}$$ is the composition
$$T_{\nu;k}=S_{\nu'}E\otimes \bigwedge^{N(\nu;k)}F\otimes S_{p(\nu;k)}G^*\to S_{\rho'}E\otimes S_{(\nu_j-k+1)'}E\otimes \bigwedge^{N(\nu;k)}F\otimes S_{p(\nu;k)}G^*
$$$$
\to S_{\rho'}E\otimes \bigwedge^{\nu_j-k+1}F\otimes \bigwedge^{N(\nu;k)}F\otimes S_{p(\nu;k)}G^*\to S_{\rho'}E\otimes \bigwedge^{N(\rho;\nu_j)}F\otimes S_{p(\nu;k)}G^*=T_{\rho;\nu_j},$$
where the first map is the Pieri map, the second  is $\bigwedge^{\nu_j-k+1}\phi$, and the third is exterior multiplication.

Finally we can also describe the terms between which we have a $(\psi\phi
)$-component map. Consider the term
$$ T_{\nu;k}=S_{\nu'}E\otimes \bigwedge^{N(\nu;k)}F\otimes S_{p(\nu;k)}G^*.$$
Consider a corner box of the partition $\nu'$ such that we can also
subtract the corresponding box from $p(\nu;k)$ in such way that we get another
nonzero term, with the same cohomology group, in the complex ${\bf
F}_\bullet$. The exterior power $\bigwedge^{N(\nu;k)} F$ will be
unaffected. The new term will occur in degree by one smaller in ${\bf
F}_\bullet$ than the original term (we decreased $\nu'$ by one box, but
the homogeneous degree from $\psi$ and the number of cohomology group
stayed the same). Between these two terms we have a $(1,1)$ degree map from
$(\psi\phi )$-component.
In other words, let $\epsilon_j$ represent the $(g-1)$-tuple with  $1$ in position $j$ and zero everywhere else. 
The map $T_{\nu;k}\to T_{\nu-\epsilon_j;k}$
is defined provided $\nu-\epsilon_j$ is a partition and $\nu_j\neq k$. The hypothesis ensures that
$$p(\nu-\epsilon_j;k)=p(\nu;k)-\epsilon_{J},$$ where 
$$J=\cases{ j
&{\rm if} $\nu_j> k$\cr
j+1&{\rm if} $k>\nu_j$.}$$
The map is the composition:
 $$ T_{\nu;k}=S_{\nu'}E\otimes \bigwedge^{N(\nu;k)}F\otimes S_{p(\nu;k)}G^*
\to S_{(\nu-\epsilon_j)'}E\otimes S_1E\otimes  \bigwedge^{N(\nu;k)}F\otimes S_1G^*\otimes S_{p(\nu-\epsilon_j;k)}G^*$$
$$ \to S_{(\nu-\epsilon_j)'}E\otimes  \bigwedge^{N(\nu;k)}F\otimes S_{p(\nu-\epsilon_j;k)}G^*=T_{\nu-\epsilon_j;k}.
$$
The first arrow is two Pieri maps
to split one box from each of  $\nu'$ and $p(\nu;k)$. The second arrow has  two components. The first component uses the map $E\otimes
G^*\rightarrow A$ given by the composition $\psi\phi$. The second component uses the maps $\psi$ and $\phi$ separately. 
To be more explicit, notice that the representation $E\otimes\bigwedge^i
F\otimes G^*$ occurs with multiplicity 2 in
$$\bigwedge^i F\otimes (E\otimes F^* )\otimes (F\otimes G^* ).$$
The two components of the second arrow involve the two possible embeddings of $\bigwedge^i F$ into
$\bigwedge^i F\otimes  F^* \otimes F$. Let us describe these two embeddings
explicitly.
We define two linear maps ${\rm tr}:{\bf K}\rightarrow F^*\otimes F$ sending $1$ to
$\sum_{i=1}^f v_i^*\otimes v_i$ for some basis $\lbrace v_1 ,\ldots
,v_f\rbrace$ of $F$.
The other is the evaluation ${\rm ev}: F\otimes F^*\rightarrow {\bf K}$. Two embeddings
of $\bigwedge^i F$ into
$\bigwedge^i F\otimes  F^* \otimes F$ are then defined as follows.
One is just
$$i_1 :=1\otimes {\rm tr} :\bigwedge^i F\to \bigwedge^i F\otimes F^*\otimes F,$$
the other is the composition
$$i_2 : \bigwedge^i F{\buildrel {\Delta\otimes {\rm tr}}\over\rightarrow}
\bigwedge^{i-1} F\otimes F\otimes F^*\otimes F\buildrel{\sigma
(2,4)}\over\rightarrow\bigwedge^{i-1}F\otimes F\otimes F^*\otimes
F{\buildrel{m\otimes 1\otimes 1}\over\rightarrow}\bigwedge^i F\otimes
F^*\otimes F$$
where $\sigma (2,4)$ switches the second and fourth factor, and $m$ denotes
the exterior multiplication.

Thus the $\phi$ and $\psi$ components of our differential are easy to
identify (up to scalar). The only problem is the $(\psi\phi )$ component
where we do not know which linear combination of maps $i_1 ,i_2$ to choose.
This problem can be solved, however, by looking at the construction of the
complex ${\bf F}_\bullet$ given in Proposition 9.

Let us choose two partitions $\lambda$ and $\nu$ such that
$\nu\subset\lambda$, $|\lambda /\nu |=1$. We have the induced map of
sheaves
$$S_{\lambda'}E\otimes {\hat{\cal M}}(\lambda )\rightarrow S_{\nu'}E\otimes
{\hat{\cal M}}(\nu )$$
which is a component of the differential in ${\cal K}_\bullet$. The induced
map of sections is the equivariant homomorphism of $A$-modules
$$f(\lambda ,\nu ):S_{\lambda'}E\otimes {\hat{ M}}(\lambda )\rightarrow
S_{\nu'}E\otimes {\hat{M}}(\nu ).$$
We know that there is a unique equivariant map
$${\hat f}(\lambda ,\nu ): S_{\lambda'}E\otimes {\bar{\bf F}}(\lambda
)_\bullet\otimes_{\bar A}A \rightarrow S_{\nu'}E\otimes {\bar{\bf F}}(\nu
)_\bullet\otimes_{\bar A}A$$
 of the minimal resolutions covering the map $f(\lambda ,\nu )$.

\proclaim Proposition 10. The $(\psi\phi )$-components of the differential
of the complex ${\bf F}_\bullet$ are the corresponding components of the
maps ${\hat f}(\lambda ,\nu )$.

\noindent{\sl Proof.} By definition of iterated cone construction. $\bullet$

\beginsection 3. The resolution of $B/J$.

Now that we have some data involving the resolution of $A/I$, we apply it
to find the terms of the resolution of $B/J$.

\proclaim Theorem 5. We have the isomorphisms ${\rm Tor}_i^A (A/I, {\bf K})={\rm Tor}_i^B
(B/J, {\bf K})$ preserving the ${\rm SL}(E)\times {\rm SL}(F)\times {\rm SL}(G)$ representation
structure
and homogeneous bidegree.

\noindent{\sl Proof.}
Consider the minimal graded free resolution of $B/J$ as an
$B$-module:
$${\bf G}_\bullet : 0\rightarrow {\bf G}_{eg+1}\rightarrow\ldots\rightarrow
{\bf G}_1\rightarrow {\bf G}_0 .$$
The complex ${\bf G}_\bullet\otimes_B B/aB$ has the $i$-th homology module
equal to ${\rm Tor}_i^B (B/J, B/aB)$. On the other hand,  
the long exact sequence of homology which is obtained by applying
$B/J\otimes_B -$ to 
the short exact sequence
$$0\rightarrow B{\buildrel a\over\rightarrow} B\rightarrow B/aB\rightarrow
0$$yields ${\rm Tor}_i^B (B/J ,B/aB)=0$ for
$i\ge 2$ and yields the exact sequence $$0\rightarrow {\rm Tor}_1^B (B/J ,B/aB)\rightarrow B/J{\buildrel
a\over\rightarrow} B/J\rightarrow B/(a,J)\rightarrow 0.$$
We know from Proposition 1 that $a$ is a nonzerodivisor on $B/J$; so, ${\rm Tor}_1^B (B/J ,B/aB)$ is also zero and  ${\bf G}_\bullet\otimes_B B/aB$ is an $A$-free resolution of
$A/I$. This resolution is minimal because the matrices of the maps in this
complex are obtained from those of maps of ${\bf G}_\bullet$ by
specializing $a$ to zero.
The terms of both minimal resolutions ${\bf G}_\bullet$ and ${\bf
G}_\bullet\otimes_B B/aB$ are the same, and they (after tensoring with ${\bf K}$)
give us the Tor groups mentioned in the theorem.$\bullet$

\proclaim Corollary. The terms in the  minimal graded free resolution, ${\bf G}_\bullet$, of the universal ring $C=B/J$ as a
$B$-module  are exactly the same as the terms of the resolution 
${\bf F}_\bullet$ of Proposition 7, once ``$A$'' is replaced by ``$B$''.

\noindent {\bf Acknowledgment.}
Thank you to Alexandre Tchernev for getting us started on this project.

\beginsection References.

\frenchspacing

\item{[AKM]} L. Avramov, A. Kustin, and M. Miller, ``Poincar\'e
series of modules over local rings of small embedding codepth or small linking
number,'' {\sl J. Alg.  \bf  118} (1988), 162--204.

\item{[B83]}
W. Bruns, ``Divisors on varieties of complexes,''
{\sl Math. Ann. \bf  264}
(1983), 53--71.

\item{[B84]}
W. Bruns,
``The existence of generic free resolutions and related objects,''
{\sl Math. Scand. \bf  55}
(1984), 33--46.

\item{[BE]}D. Buchsbaum and D. Eisenbud, ``Some structure theorems for finite free resolutions'' {\sl Advances Math. \bf 12} (1974)  84--139.

\item{[E]}D. Eisenbud
``Commutative Algebra with a view toward Algebraic Geometry,''  
Spring\-er Verlag, 
 New York (1995).

\item{[Ha]}  M. Hashimoto, \hskip-.61pt ``Determinantal ideals without
minimal free resolutions,''  \hskip-.61pt {\sl Nagoya Math. J. \bf  118} (1990), 
203--216.

\item{[He]} R. Heitmann, 
``A counterexample to the rigidity conjecture for rings,'' {\sl Bull. Amer. Math. Soc. (N.S.) \bf  29
} (1993), 94--97.

\item{[Ho]}
M. Hochster,  ``Topics in the homological theory of modules over commutative rings'', 
CBMS Regional Conf. Ser. in Math., no. 24
(1975), 
Amer. Math. Soc., Providence, RI. 

\item{[Hu]} C. Huneke, 
``The arithmetic perfection of Buchsbaum-Eisenbud varieties and generic modules of projective dimension two,'' {\sl
Trans. Amer. 
Math. Soc. \bf  265 } (1981), 211--233.

\item{[K]} A. Kustin,  ``The resolution of the universal ring for finite length modules of projective dimension two'', preprint, http://www.math.sc.edu/$\sim$kustin/research.html.

\item{[KM]} A. Kustin and M. Miller, ``Multiplicative
structure on resolutions of algebras defined by Herzog ideals,'' {\sl J. London
Math. Soc. (2) \bf  28 
} (1983), 247--260.

\item{[PW]}P. Pragacz and J. Weyman, ``On the generic free
resolutions,'' {\sl J. Alg. \bf  128} (1990), 1--44.

\item{[R]} P. Roberts, ``Homological Invariants of Modules over Commutative Rings,'' 
 Les Presses de l'Universit\'e de Montr\'eal,
Montr\'eal,
(1980).

\item{[T]} 
A. Tchernev, ``Universal complexes and the generic structure of free resolutions,'' {\sl Mich. Math. J. \bf  49} (2001), 65--96.

\item{[W]} 
J. Weyman, ``Cohomology of vector bundles and syzygies,'' Cambridge University Press,  Cambridge (2003).

\bye